\documentclass[12pt]{article}
\usepackage{latexsym}
\usepackage{amsfonts}
\usepackage{amsmath}
\usepackage{amsthm}
\usepackage{epsfig}

\newtheorem{question}{Question}

\newtheorem{thm}{Theorem}
\newtheorem{defn}[thm]{Definition}
\newtheorem{lemma}{Lemma}

\newcommand{\y}{\bar{y}}
\newcommand{\Z}{\mathbb Z}
\newcommand{\R}{\mathbb R}
\newcommand{\N}{\mathbb N}
\newcommand{\E}{{\bf E}}
\newcommand{\prob}{{\bf P}}
\newcommand{\eps}{\epsilon}
\newcommand{\be}{\begin{equation}}
\newcommand{\ee}{\end{equation}}

\newcommand{\var} {{\mbox{var}}}
\newcommand {\spec} {specification }

\newcommand {\specs} {specifications }
\newcommand {\kappak} {K}
\newcommand {\fernandez} {$\text{Fern}\acute{\text{a}}\text{ndez }$}

\def\23#1{\lfloor 2{#1}/3 \rfloor}

\def\22#1{2^{2^{#1}}}

\title{Nonuniqueness for specifications in $\ell^{2+\epsilon}$}
\author{Noam Berger \and Christopher Hoffman\thanks{supported in part by NSF grant \#
62-5254} \and Vladas Sidoravicius }

\begin{document}
\maketitle
\begin{abstract}
For every $p>2$, we construct a regular and continuous
specification ($g$-function), which has a  variation sequence that
is in $\ell^p$ and which admits multiple Gibbs measures. Combined
with a result of Johansson and \"Oberg, \cite{sweden}, this
determines the optimal modulus of continuity for a specification
which admits multiple Gibbs measures.

\end{abstract}

\section{Introduction}\label{sec:intro}
For a finite set $A$ let $P(A)$ be the set of probability
distributions on $A$. A {\it specification g} (also commonly known
as a $g$-function) is a measurable function from $A^\N$ to $P(A)$.
A specification $g$ is {\it regular} if there exists $\varrho>0$
such that for every sequence $\{b_n\}_{n\in\N}$ and every $a\in A$
we have that $(g(\{b_n\}))(a)\geq\varrho$.  We focus on
specifications that are regular and continuous with respect to the
product topology.
%
\noindent A {\it Gibbs measure} for a specification $g$ is a shift
invariant probability measure $\mu$ on $A^\Z$ such that
 for every $f:A\to \R$,
\begin{equation*}
\E_\mu\left( f(x_0)|x_{-1},x_{-2},\ldots \right) =
\E_{g(x_{-1},x_{-2},\ldots)}(f) \ \ \ \ \ \ \ \mbox{a.s.}
\end{equation*}
It is easy to show that every continuous specification has a Gibbs
measure.

Given a past, $x_{-1},x_{-2},\dots$, the specification $g$ 
tells us the
probability distribution for $x_{0}$, the next state of the process.
Thus the specification and the past determine the stochastic evolution
of the process.  One common example of a specification is a finite
state space Markov chain.  The specification for a $k$-step Markov
chain is determined by $x_{-1},\dots,x_{-k}$.  For this reason
D\"oblin and Fortet referred to specifications as ``chains
with infinite connections'' \cite{DF}.

The question of whether a specification
$g$ uniquely  determines the stationary process under a natural
``mixing'' assumption has been prominent 
since the pioneering work of D\"oblin and Fortet \cite{DF}. 
However, in the last
three decades reasonable progress was achieved only in providing
sufficient conditions for the uniqueness of Gibbs measure (see
\cite{stenflo2} for detailed discussions and references). Harris
\cite{H} studied the behavior of lumped Markov chains and
introduced important coupling ideas that were used by several
later authors. M. Keane \cite{keane} introduced the notion of a
continuous $g$-function and gave conditions under which a
$g$-function has a unique measure. One natural way to express
uniqueness conditions is in terms of the modulus of continuity of
$g$.
To quantify the modulus of continuity of a specification $g$, we
define the {\it variation of $g$ at distance $k$} to be
\[
\var_k(g)=\sup\left\{ \|g(b)-g(b')\|_1 \ \left| \
b_1=b'_1,b_2=b'_2,\ldots,b_k=b'_k \right.\right\}.
\]
The continuity of $g$ is equivalent to the sequence $\var_k(g)\to
0$.

In \cite{walters}, Walters showed that if $\var(g)\in\ell^1$, then
the Gibbs measure is unique. Walters and Ledrappier
\cite{francois} established a strong connection between \specs and
the thermodynamic formalism of statistical mechanics.  In
particular they relate \specs with symbolic dynamics as developed
in works of Sinai \cite{yag}, Ruelle \cite{d1}, Bowen \cite{Bo},
and others.

Walters's work was sharpened by Lalley \cite{lalley} and Berbee
\cite{berbee}, who showed that the Gibbs measure is unique if
\[
\sum_{n=1}^\infty \exp\left(-\sum_{m=1}^n\var_m(g)\right)=\infty.
\]

On the other hand, by the early 1980's it was known that the
equilibrium measures of appropriately chosen one dimensional
long-range Ising models are not unique (see \cite{Dy} for the
hierarchical type models and \cite{FS} for the $1/r^2$ decaying
cases). Existence of phase transition was later established for
one-dimensional long-range percolation model (see \cite{NS} and
\cite{AN}) and finally for one-dimensional FK-Random Cluster model
\cite{ACCN}. Nevertheless a little more than a decade ago it was
widely believed that continuous and regular specifications admit a
unique Gibbs measure.

However, in 1993, Bramson and Kalikow \cite{bk} provided a
remarkable (and until now unique) example of a continuous and
regular specification that admits multiple Gibbs measures.
The variation of the function $g$ in Bramson and Kalikow's
construction is not in $\ell^p$ for any $p$. In fact, in their
example $\var_k(g)\geq\frac{C}{\log k}$ for some constant $C$.
This gave rise to the following question: For which values of
$p\,$ does $\var(p)\in\ell^p$ imply uniqueness.

A few years ago Stenflo further sharpened Berbee's work
\cite{stenflo}. However, the results of Berbee and Stenflo, while
improving over Walters' result, are still in the realm of
$\ell^1$.  Recently Johansson and \"Oberg \cite{sweden} showed
that if $g$ is regular and $\{\var_k(g)\}_{k=1}^\infty$ is in
$\ell^2$ then $g$ admits a unique Gibbs measure. Our main result
is the following:

\begin{thm}\label{thm:main}
For every $p>2$, there exists a regular specification $g$ such
that $\{\var_k(g)\}\in\ell^p$ and $g$ admits multiple Gibbs
measures.
\end{thm}
This shows that the result of \"Oberg and Johansson is tight.

\noindent {\bf Remark.} In \cite{fern1,fern2} \fernandez and
Maillard proved a Dobrushin type uniqueness condition. This
condition is  not comparable with the variation conditions.

\section{Construction}\label{sec:construct}



We will use the alphabet of size four $A=\{+1,-1\}^2$. We fix a
parameter $\epsilon \in (0,.5)$.  Given this we pick a positive
integer $\kappak=\kappak(\epsilon)$ such that the inequalities in
lines (\ref{star1'}) through (\ref{star3'}) below are true for all
$k\geq K$.  We now begin to define a regular continuous \spec
$g(x,y)=g^{\epsilon}(x,y):A^\N\to P(A).$

The choice of $(x_0,y_0)$ given $\{x_{-i},y_{-i}\}_{i=1}^\infty$
consists of four steps:
\begin{enumerate}
\item\label{item:gety}
Choose $y_0$ independently of $\{x_{-i},y_{-i}\}_{i=1}^\infty$, so
that $y_0=+1$ with probability $0.5$ and $y_0=-1$ with probability
$0.5$.
\item\label{item:getS}
Using the values of $\{y_{-i}\}_{i=0}^\infty$ choose a
(deterministic) set of odd size $S\subset -\N$.
\item\label{item:getupsilon}
Using the values of $\{y_{-i}\}_{i=0}^\infty$ choose a
(deterministic) value $0\leq\upsilon < 0.4$.
\item\label{item:getx}
Let $z$ be the majority value of $\{x_t:t\in S\}$. Choose $x_0=z$
with probability $0.5+\upsilon$ and $x_0=-z$ with probability
$0.5-\upsilon$.
\end{enumerate}

In order to complete the first step the second coordinate ($y$)
must be i.i.d.\ with distribution (1/2,1/2). We ensure this if for
all $x$ and $y$
\[
g(x,y)(1,1)+g(x,y)(-1,1)=g(x,y)(1,-1)+g(x,y)(-1,-1)=\frac12.
\]
First we pick $y_0$.  We write $\y$ to represent $y$ and $y_0$.
The most intricate part of the construction is the choice of the
set $S$.  Before we choose $S$ we need some notation.

For every positive integer $k$, let $I_k$ be the sequence of
length
\[
\ell_k=\left\lceil(1+\eps)^k\right\rceil
\]
such that the last element is $-1$ and all other elements are
$+1$. We also define
\begin{equation}\label{eq:defbeta}
\beta_k=2^{\ell_k}
\end{equation}
and
\begin{equation}\label{eq:defnu}
\nu_k=\frac{\beta_k}{\beta_{k-1}}\sim 2^{\eps\cdot\ell_{k-1}}
\end{equation}


Now we define a block structure which will allow us to choose the
 set $S\in-\N$ and a
parameter $\upsilon$.
\begin{defn}\label{compblock}
A {\bf complete $k$ block in $\y$} is a subsequence
$\{y_i\}_{i=a}^{b-1}$ such that
\begin{enumerate}
    \item $\left(\y_{a-\ell_k+1},\y_{a-\ell_k+2},\ldots,\y_a\right)=
     \left(\y_{b-\ell_k+1},\y_{b-\ell_k+2},\ldots,\y_b\right)=I_k$
    \item for no $c \in (a,b)$ does
    $\left(\y_{c-\ell_k+1},\y_{c-\ell_k+2},\ldots,\y_c\right)=I_k$
    \end{enumerate}
\end{defn}

We also define a partial $k$ block as follows.
\begin{defn}
A {\bf partial $k$ block in $\y$} is a subsequence
$\{y_i\}_{i=a}^{0}$ such that
\begin{enumerate}
    \item $\left(\y_{a-\ell_k+1},\y_{a-\ell_k+2},\ldots,\y_a\right)= I_k$
    \item for no $c>a$ does $\left(\y_{c-\ell_k+1},\y_{c-\ell_k+2},\ldots,\y_c\right)=I_k$
\end{enumerate}
\end{defn}

\begin{defn}
A {\bf $k$ block in $\y$} is an interval that is either a complete
$k$ block
or a partial $k$ block.
\end{defn}

If $\y$ is well understood then we write $B=[a,b)$ to denote a $k$
block $\{y_i\}_{i=a}^{b-1}$. Note that this definition is
invariant under a shift of $\y$ in the following sense: Given $y$
and any integer $t<0$ define $\sigma^{t}(\y)$ by
$\sigma^{t}(\y)_i=y_{i+t}$ for all $i\leq 0$. Then for $b<0$, if
$[a,b)$ is a $k$ block for $\sigma^t(\y)$ then $[a+t,b+t)$ is a
$k$ block for $\y$. Also note that the length of a $k$ block is
likely to be close to $\beta_k$, and that the number of $k-1$
blocks inside a $k$ block is likely to be close to $\nu_k$. A
precise statement of these claims will be used extensively in the
next section.

Next we label all of the $k$ blocks.  Given a sequence $\y$ and a
positive integer $k$ we will define $B_{k,i}=B_{k,i}(\y)$ to be
the $i$-th $k$ block in $\y$. More precisely we define
$a_{k,i}=a_{k,i}(\y)$ and $b_{k,i}=b_{k,i}(\y)$ such that
\begin{enumerate}
    \item $[a_{k,i},b_{k,i})$ is a $k$ block in $\y$.
    \item $b_{k,i+1} = a_{k,i}$ for all $i$
    \item $b_{k,1}=1$.
\end{enumerate}
These sequences can be either finite or infinite. Then we define
$B_{k,i} = [a_{k,i},b_{k,i}).$

Given $y$ and $k > \kappak$ define $N_k(\y)$ to be the number of
$k-1$ blocks in the $k$ block containing $0$. More precisely, we
take $N_k(\y)$ to be so that
\begin{equation}\label{defnk}
a_{k,1}=a_{k-1,N_k(\y)},
\end{equation}
or $N_k=\infty$ if (\ref{defnk}) has no solution. We also set
$N_{\kappak}(\y)=|a_{\kappak,1}(\y)|.$



\begin{defn}\label{def:begin}
We say that the {\bf beginning} of a $\kappak$ block $B=[a,b)$ is
the interval
$O(B)=\left[a,\min(b,a+\beta_\kappak^{1-\epsilon})\right)$.

For $k > \kappak$ we say that the {\bf beginning} of a $k$ block
$B=[a,b)$ is the $\nu_k^{1-\epsilon}$ first $k-1$ blocks in $B$,
i.e.
\begin{equation*}
O(B) = \left\{
\begin{array}{ll}
[a,b) & \mbox{ if } N_k(\sigma^b(\y))\leq\nu_k^{1-\epsilon} \\
{}[a,b_{k-1,N_k(\sigma^b(\y)) - \nu_k^{1-\epsilon}}(\sigma^b(\y)))
& \mbox{otherwise }
\end{array}
\right..
\end{equation*}
\end{defn}
\begin{defn}
The {\bf opening} $C(B)$ of a $k$ block $B$ is the set of points
$t\in B$ such that
\begin{enumerate}
\item
$t$ is in the beginning of its $j$-block for every $\kappak\leq
j\leq k$.
\item
If $a_j(t)$ is the smallest element in the $j$-block containing
$t$, then $t-a_j(t)<\beta_{j+1}$ for all $\kappak\leq j\leq k$.
\end{enumerate}
\end{defn}

\begin{figure}[h]
\center{\epsfig{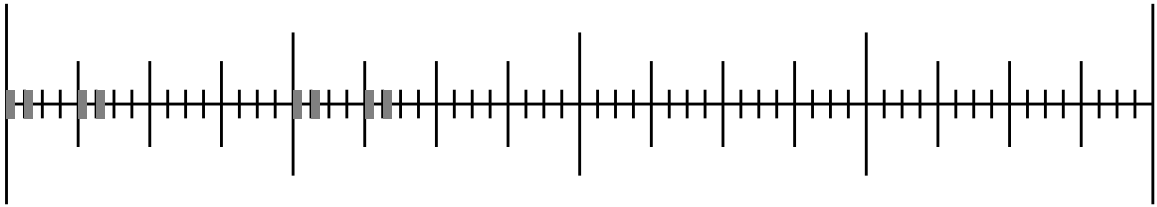}} \caption{A block.
The opening is marked in gray.}\label{fig:open}
\end{figure}

We also define $C_k$ to be the union of $C(B)$ over all $k$ blocks
$B$. Note that if $t\in C_k$ then $t\in C_j$ for all $j<k$. The
event that $0 \in C_k$ is determined by $\y_i$, $i \in
[-\beta_{k+1},0]$. This fact will be used to show that our \spec
is continuous and to show how quickly $\text{var}_k(g)$ approaches
0. Define $k_0$ to be the highest value such that $0\in C_{k_0}$.
If there is no such value then we take $k_0=\kappak-1$. If $0\in
C_k$ for every $k$ then we take $k_0 = \infty$.

We now define $S$ to be the following set: If $k_0=\infty$ then
$S$ is the empty set.
If $|C(B_{k_0+1,1})|$ is odd then we take $S=C(B_{k_0+1,1})$. If
$|C(B_{k_0+1,1})|$ is even,
 $S=C(B_{k_0+1,1}) \setminus \{\max(C(B_{k_0+1,1}))\}$.

Next we choose the value of $\upsilon$. If
$\min(S)<-\beta_{k_0+2}$ then we take $\upsilon = 0$. Otherwise,
we take
\begin{equation}\label{eq:defups}
\upsilon=\left\{
\begin{array}{ll}
0.4  & \mbox{ if }k_0=\kappak-1\\
\beta_{k_0}^{-\frac{1}{2}+\epsilon} & \mbox{ if } k_0\geq\kappak.
\end{array}
\right.
\end{equation}
Finally we set
\begin{equation}
g_1(x,y,y_0)=0.5+\upsilon\cdot\mbox{sign}\sum_{t\in S}x_t
\end{equation}
where sign($x$) is the function that takes on values +1 or -1
depending on whether $x$ is positive or negative.
%
Thus we have defined our \spec
$g$ by the coordinate functions
\begin{equation*}
\begin{array}{rll}

g_{(1,1)}(x,y)   & = &   .5     g_1(x,y,1)  ,                      \\
g_{(-1,1)}(x,y)  & = &   .5 (1- g_1(x,y,1)) ,                      \\
g_{(1,-1)}(x,y)  & = &   .5     g_1(x,y,-1) ,                      \\
g_{(-1,-1)}(x,y) & = &   .5 (1- g_1(x,y,-1)),                      \\
\end{array}%
\end{equation*}

We leave it to the reader to check the following lemma which says
that $g$ is symmetric in the $x$ coordinate.
\begin{lemma}
\label{lem:symmetry} For all $x,y \in \{1,-1 \}^{-\N}$ and $a,b
\in \{1,-1\}$
$$ g_{(a,b)}(x,y)=g_{(-a,b)}(-x,y).$$
\end{lemma}

\subsection{Continuity}

We now show that $g$ is continuous.  The two most important
elements of the construction are that for all $y$
\begin{enumerate}
\item there exists $N$ such that the set $S$ is determined by $\{y_i\}_{i=-N}^{0}$.
\item The larger $N$ is the closer $\upsilon$ is to zero.
\end{enumerate}
In the next lemma we quantify these two statements and show that
they imply that $g$ is continuous.

\begin{lemma}\label{lem:cont}
The \spec
$g$ is regular and continuous. Moreover the sequence
$\text{var}_j(g) \in l^p$ for all
\begin{equation}\label{eq:noeps}
p>\frac{2(1+\epsilon)^2}{1-2\epsilon}.
\end{equation}
\end{lemma}
Note that as $\epsilon$ approaches $0$, the bound in
(\ref{eq:noeps}) goes to $2$.
\begin{proof}
%
%
%
Let $k > \kappak$. We want to estimate $\var_j(g)$ for
$\beta_k<j\leq \beta_{k+1}$.

Let $\{x^{(1)}_i,y^{(1)}_i\}_{i\in-\N}$ and
$\{x^{(2)}_i,y^{(2)}_i\}_{i\in-\N}$ be such that
$x^{(1)}_i=x^{(2)}_i$ and $y^{(1)}_i=y^{(2)}_i$ for every $i>-j$.
It is enough to estimate
$|g_1(x^{(1)},y^{(1)},l)-g_1(x^{(2)},y^{(2)},l)|$ for $l=-1,+1$.
Fix $l$ and for $h=1,2$ let
\[
\y^{(h)}_i=\left\{
\begin{array}{ll}
y^{(h)}_i & \mbox{ if } i < 0 \\
l & \mbox{ if } i = 0
\end{array}
\right..
\]
If $k_0(\y^{(1)}) \leq k-2$ then $k_0(\y^{(1)})=k_0(\y^{(2)})$ and
$g_1$ depends only on $l,\{x^{(h)}_i,y^{(h)}_i\}_{i>-\beta_k}$.
Therefore, on this case
$g_1(x^{(1)},y^{(1)},l)-g_1(x^{(2)},y^{(2)},l)=0$.

If, on the other hand, $k_0(\y^{(1)}) \geq k-1$, then
$k_0(\y^{(2)}) \geq k-1$ as well and therefore
$|g_1(x^{(h)},y^{(h)},l)-0.5|<\beta_{k-1}^{-\frac12+\epsilon}$.
Thus
\begin{equation}
\text{var}_j(g)<
 2\left(\beta_{k-1}\right)^{-\frac12+\epsilon}<
 8\left((\beta_{k+1})^{1/(1+\epsilon)^2}\right)^{-\frac12+\epsilon}<
 8j^{-\frac{1-2\epsilon}{2(1+\epsilon)^2}}.
\end{equation}
\end{proof}

\section{Multiple measures}\label{sec:mult}

The goal of this  section is to prove the following lemma:
\begin{lemma}\label{prop:main}
For every $\epsilon \in (0,.5)$  the \spec $g^{\epsilon}$ admits
multiple measures.
\end{lemma}
Theorem \ref{thm:main} is a direct consequence of Lemmas
 \ref{lem:cont} and \ref{prop:main}.


To see that the function $g$ admits multiple Gibbs measures, we
introduce the following notation. Choose an arbitrary $t$. We let
$B_k(t)$ denote the $k$ block containing $t$, and we let
$C(B_k(t))$ denote the opening of
$B_k(t)$. 
%
%
We will prove that $X(B_k(t))=X(B_{k+1}(t))$ with extremely high
probability.
Next we note that for any $t$ and $t'$ that $B_k(t)=B_k(t')$ (and
thus $X(B_k(t))=X(B_k(t'))$) for all $k$ large enough. (See Figure
\ref{fig:arrows} below.)  Thus by the symmetry of Lemma
\ref{lem:symmetry} there exist at least two invariant measures:
one where $\lim_{k\to\infty}X(B_k(t))=1$ a.s.\ and one where
$\lim_{k\to\infty}X(B_k(t))=-1$ a.s.


\begin{figure}[h]
\center{\epsfig{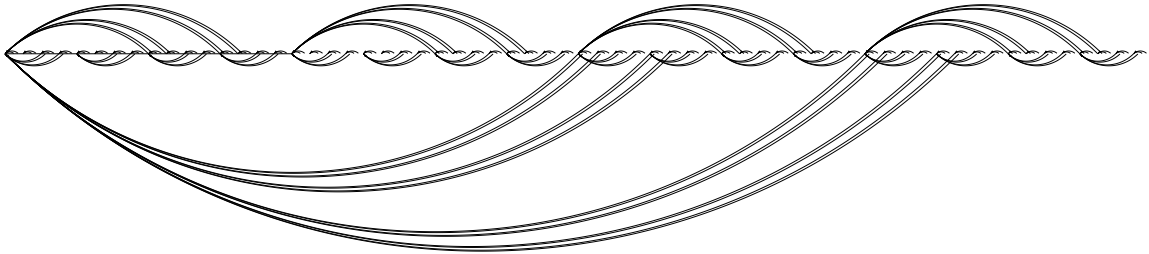}} \caption{The arrows
point from every point $t$ to a point in
$S(t)$.}\label{fig:arrows}
\end{figure}

Assume that $\mu$ is a stationary and ergodic measure on
$\{(\pm1,\pm1)\}^{\Z}$ which is a Gibbs measure for $g$.
Note that in the definitions of a complete $k$ block
(Definition \ref{compblock}
of Page \pageref{compblock}) there was no requirement that $a$ and
$b$ be negative. Thus we can use the same definition of a complete
block for a sequence in $\{\pm 1 \}^{\Z}$

\noindent {\bf Remark:}
From now on we will abbreviate the term {\bf complete $k$ block}
to {\bf $k$ block}. This is a natural abuse of notation that
results from the fact that we now speak about objects taking
places in $\{(\pm1,\pm1)\}^{\Z}$ instead of
$\{(\pm1,\pm1)\}^{-\N}$. In the same spirit, $x$ and $y$ will now
denote two-sided sequences.

In subsections \ref{subsec:gb} and \ref{subseq:beautipoints} we
prove results about the block structure. The block structure
depends only on the $y$ sequence, which is an i.i.d.\ $(0.5,0.5)$
sequence. In subsection \ref{subsec:last} we use the results
obtained in subsections \ref{subsec:gb} and
\ref{subseq:beautipoints} to understand the (more complicated)
structure of the $x$ sequence and show the existence of two
different Gibbs measures.

\subsection{Good blocks}\label{subsec:gb}
In this subsection we define the notion of a good block. A good
block is a block whose length is close to its expected length and
the size of its opening is also close to its expected value.
Our aim in this subsection is to show that a $k$ block is good
with high probability. In the next subsections we show that good
blocks exhibit behavior that yields multiple Gibbs measures.
\begin{defn}\label{defn:goodblock}
A $k$ block $B$ is {\bf good} if
\begin{enumerate}
\item\label{item:size}
$|B| < \beta_k^{1+\epsilon}$ and
\item\label{item:coresize}
$|C(B)|>\beta_k^{1-\epsilon}2^{-k}$
\end{enumerate}
\end{defn}

In order to show that most blocks are good we put a measure on $k$
blocks by conditioning on $y$ such that
$\{y_i\}_{i=-l_k+2}^{1}=I_k$.  For a $k$ block
$B=\{y_i\}_{i=1}^{n-1}$ such that
\begin{enumerate}
\item $y_1=-1$,
\item $y_{n-l_k}=\dots y_{n-1}=1$ and
\item there is no occurrence of $I_k$ in $B$
\end{enumerate}
we define the measure
$$m_k(B)=\prob(B_k(1)=B|\{y_i\}_{i=-l_k+2}^{1}=I_k)=2^{-|B|}.$$
For $t\in\Z$, we let $B_k(t)$ denote the $k$ block containing $t$.
With these definitions we can check that for any $\mu$ which is a
Gibbs measure for $g$ and any $k$ and $B$
$$\prob(B_k(0) \mbox{ is a translate of $B$})=\frac{|B|m_k(B)}{\beta_k}.$$

\begin{lemma} \label{lem:pvg}
For all $k \geq K$ $$\sum_{B \mbox{ not
good}}m_k(B)<3\cdot2^{-3k}.$$
\end{lemma}

\begin{proof}
We prove the lemma by induction on $k$. For $k=K$ we need only to
check the first condition.  The probability that a $k$ block is
longer than $c\ell_k$ is less than $(1-2^{-\ell_k})^c$. Thus the
probability that the first condition in the definition is not
satisfied is at most
\begin{equation}\label{star1'}
(1-2^{-\ell_k})^{\beta_k^{1+\epsilon}/2\ell_{k-1}}<2^{-3k}.
\end{equation}
This is true by our choice of $K$.

If $k>K$ and
\begin{enumerate}
\item $|B| < \beta_k^{1+\epsilon}$
\item there are at least $\nu_k^{1-\epsilon}$ k-1 blocks in $B$
and \label{ii}
\item at least half of the $\nu_k^{1-\epsilon}$ k-1 blocks in the beginning of
$B$ are good \label{iii}
\end{enumerate}
then
\begin{equation} \label{eq:star1}
|C(B)|
\geq\frac{1}{2}\nu_k^{1-\epsilon}\beta_{k-1}^{{1-\epsilon}}2^{-k+1}
\geq\beta_k^{{1-\epsilon}}2^{-k}
\end{equation}
and $B$ is good.

As the bound in Line \ref{star1'} holds for any $k>K$ we have that
the probability that the first condition is not satisfied is at
most $2^{-3k}$. The probability that conditions \ref{ii} and
\ref{iii} are satisfied is the probability that if we select
$\nu_k^{1-\epsilon}$ $k-1$ blocks independently according to $m$
that none of them are a $k$ block and that at least half of them
are good. The probability that a $k-1$ block is also a $k$ block
is $\nu_k^{-1}$. Thus the probability of having a $k$ block among
$\nu_k^{1-\epsilon}$ k-1 blocks is less than
\begin{equation} \label{eq:star2}
\nu_k^{1-\epsilon}\nu_k^{-1}<\nu_k^{-\epsilon}\leq 2^{-3k}.
\end{equation}

By the induction hypothesis a $k-1$ block is not good with
probability less than $2^{-k+1}$.  Thus the probability that half
of a sequence of $\nu_k^{1-\epsilon}$ $k-1$ blocks chosen
independently are not good is at most

\begin{equation} \label{star3'}
2^{\nu_k^{{1-\epsilon}}}(2^{-k})^{.5\nu_k^{1-\epsilon}} \leq
2^{-\nu_k^{{1-\epsilon}}}\leq 2^{-3k} .
\end{equation}
\end{proof}


\begin{lemma}\label{lem:goodhighp}
 For every  $k\geq K$
 $$\prob(B_k(0) \mbox{ is good})>1-2^{-k}.$$
\end{lemma}

\begin{proof}
First, we want to estimate
\[
\tilde{Y}=\sum_{|B|\leq\beta_k}m_k(B).
\]
$\tilde{Y}$ is the probability that $I_k$ appears somewhere in
$\{y_j\}_{j=1}^{\beta_k}$. For $i=1,\ldots,\ell_k$, let $E_i$ be
the event that there exists $h$ with $h\equiv i \mod \ell_k$ such
that $(y_h,y_{h+1},\ldots,y_{h+\ell_k-1})=I_k$. $\{E_i\}$ are
negatively associated, and for every $i$,
\[
\prob(E_i)=1-(1-\beta_k^{-1})^{\beta_k/\ell_k}.
\]
Therefore
\begin{equation}\label{eq:geqyt}
\tilde{Y}\geq 1-(1-\beta_k^{-1})^{\beta_k}.
\end{equation}
On the other hand, for $i=1,\ldots,1000$ let $F_i$ to be the event
that $I_k$ appears in
\[
\{y_j:j=\beta_ki/1000\ldots\beta_k(i+1)/1000-1\}.
\]
Then the $F_i$-s are independent and using a first-moment argument
we see that for every $i$, $\prob(F_i)\leq 1/1000$. The event
measured by $\tilde{Y}$ is the union of the $F_i$-s plus the event
that $I_k$ appears on the seams between the blocks. Therefore,
\begin{equation}\label{eq:leqyt}
\tilde{Y}\leq 1-(1-1/1000)^{1000}+1000\ell_k2^{-\ell_k}.
\end{equation}

For any $j$, let
\[
\tilde{Z}=\sum_{j\beta_k<|B|\leq(j+1)\beta_k}m_k(B)
\]
be the probability that the first appearance of $I_k$ is in the
interval $[j\beta_k,(j+1)\beta_k)$. Then
\begin{equation}\label{eq:conyz}
\tilde{Z}=\tilde{Y} \left( 1-\tilde{Y} \right)^j \left(
1-\ell_k2^{-\ell_k} \right)^j,
\end{equation}
where the last term comes from the event that $I_k$ appears on the
seam between two consecutive intervals of length $\beta_k$.

From (\ref{eq:geqyt}), (\ref{eq:leqyt}) and (\ref{eq:conyz}), we
get that for any $j$ we have that
$$(5/2)^{-j}> \sum_{j\beta_k<|B|\leq(j+1)\beta_k}m_k(B)>2^{-2j-2}.$$
From this we get
$$\sum_{|B|>k\beta_k}m_k(B)>\sum_{j\geq k}2^{-2j-2} >3\cdot2^{-3k}$$
and
$$\sum_{|B|>k\beta_k}|B|m_k(B)<\sum_{j\geq k}(5/2)^{-j}(j+1) \beta_k < 2^{-k}\beta_k.$$
Thus for any set $S$ of $k$ blocks such that
$$\sum_{B\in S}m_k(B)<3\cdot 2^{-3k}$$
we have that
$$\sum_{B\in S}|B|m_k(B)<2^{-k}\beta_k.$$
Combining this last statement with Lemma \ref{lem:pvg}
 $$\prob(B_k(0) \mbox{ is good})\geq
 1-\frac1{\beta_k}\left(\sum_{B\mbox{ not good}}|B|m_k(B)\right)
 >1-2^{-k}.$$
\end{proof}

\subsection{Beautiful points}\label{subseq:beautipoints}
In this subsection we define the notion of a beautiful point. Our
goal in this subsection will be to show that most points are
beautiful.

\begin{defn}\label{def:beau}
We say that $t\in\Z$ is {\bf $k$ beautiful} if for every $j\geq
k$,
\begin{enumerate}
\item\label{item:ingood}
$B_j(t)$ is good and
\item\label{item:inend}
$t$ is not in the beginning of $B_j(t)$. (The {\em beginning} of a
block was defined in Definition \ref{def:begin} of page
\pageref{def:begin})
\end{enumerate}
\end{defn}

In this subsection we state two easy lemmas:
\begin{lemma}\label{lem:triv}
If $t$ and $s$ belong to the same $k$ block, and $t$ is $k+1$
beautiful, then $s$ is also $k+1$ beautiful.
\end{lemma}

\begin{proof}
Lemma \ref{lem:triv} follows immediately from Definition
\ref{def:beau}. \end{proof}

\begin{lemma}\label{lem:singlepoint}
Almost surely, for every $t$ there exists $\hat{k}(t)$ such that
$t$ is $\hat{k}(t)$ beautiful.
\end{lemma}

\begin{proof}
First we show that
$$\prob(0 \mbox{ is $k$ beautiful})>1-2^{-k+2}.$$
%
%



Lemma \ref{lem:goodhighp} tells us that the probability that the
$j$ block containing $0$ is good is greater than $1-2^{-j}$. As
there are at most $\nu_j^{1-\epsilon}$ $k-1$ blocks in the
beginning of $B_k(0)$ and the expected number of $k-1$ blocks in
$B_k(0)$ is $\nu_j$, the probability that $0$ is in the beginning
of $B_j(0)$ is less than
$$\nu_j^{1-\epsilon}/\nu_j=\nu_j^{-\epsilon}.$$
Therefore $0$ is $k$ beautiful with probability at least



\[
1-\sum_{j=k}^\infty (2^{-j}+\nu_j^{-\epsilon})=1-\sum_{j=k}^\infty
2^{-j+1}=1-2^{-k+2}.
\]

%
%

Thus by Borel-Cantelli we get that $\hat{k}(0)$ exists a.s. The
lemma is true because $\mu$ is shift invariant.
\end{proof}

\subsection{Proof of lemma \ref{prop:main}}\label{subsec:last}
In the previous subsections we only discussed the structure
induced by the $y$ values. In this subsection we will shed some
light on the structure of the $x$ values. For a $k$ block $B$, we
define the {\bf signature} of $B$ to be
\begin{equation}\label{eq:sign}
X(B)=\mbox{sign}\sum_{t\in C^\prime(B)}x_t
\end{equation}
where
\[
C^\prime(B)=\left\{
\begin{array}{ll}
C(B) & \mbox{ if } |C(B)| \mbox{ is odd}\\
C(B)\setminus\{\max(C(B))\} & \mbox{otherwise}
\end{array}
\right..
\]

Note that since $C^\prime(B)$ is always odd, $X(B)$ can only be
$+1$ or $-1$. We assume that $\mu$ is an ergodic Gibbs measure for
$g$.

\begin{lemma}\label{lem:majinblock}
For all $t$
$$X(t)=\lim_{k \to \infty}X(B_k(t))$$
exists a.s.\ and is equal to 1 or $-1$.
\end{lemma}
\begin{proof}
By Lemma \ref{lem:singlepoint} we get $\hat{k}$ such that $t$ is
$\hat{k}$ beautiful. By Lemma \ref{lem:symmetry} it causes no loss
of generality to assume that $X(B_{k+1}(t))=1$. Using that
$B_k(t)=[a,b)$ is good for every $k>\hat{k}$, that $\mu$ is a
Gibbs measure with respect to $g$, and line (\ref{eq:defups}) in
the definition of $g$, we get that given $y$ and $\{x_i\}_{i<a}$
the values of $x_s:s\in C^\prime(B_k(t))$ are independent and
identically distributed. For every $s$ in $C^\prime(B_k(t))$,
\begin{equation}\label{eq:sboded}
\prob(x_s=1\ |\ y, \{x_i\}_{i<a})
=\frac{1}{2}+\beta_k^{-\frac12+\epsilon}.
\end{equation}
Also, the fact that $B_k(t)$ is good implies
$$|C^\prime(B_k(t))|\geq\beta_k^{{1-\epsilon}}2^{-k}-1.$$

Thus we get
$$\prob(X(B_{k}(t)=X(B_{k+1}(t))= \prob\left(\sum_{s \in
 C^\prime(B_k(t))}x_s>0\ |\ y, \{x_i\}_{i<a}\right).$$
We have that
$$\E\left(\sum_{s \in C^\prime(B_k(t))}x_s \ |\ y, \{x_i\}_{i<a}\right)
    =2|C^\prime(B_k(t))|\beta_k^{-\frac12+\epsilon}.$$
The standard deviation of the sum is less than
$\frac12\sqrt{|C^\prime(B_k(t))|}$. Thus by Markov's inequality we
have that

\begin{eqnarray*}
\prob(X(B_{k}(t)=X(B_{k+1}(t))
 &\leq & \prob\left(\sum_{s \in C^\prime(B_k(t))}x_s>0\ |\ y, \{x_i\}_{i<a}\right)\\
 &\leq & \left(\frac{2|C^\prime(B_k(t))|\beta_k^{-\frac12+\epsilon}}
        {\frac12\sqrt{|C^\prime(B_k(t))|}} \right)^{-2}\\
 &\leq & 16|C^\prime(B_k(t))|^{-1}\beta_k^{1-2\epsilon}\\
 &\leq & 32\beta_k^{1-\epsilon}2^{k}\beta_k^{1-2\epsilon}\\
 &\leq & 32\beta_k^{-\epsilon}2^{k}\\
 &\leq & 2^{-k}.
\end{eqnarray*}

By Borel-Cantelli there are only finitely many values of $k$ such
that $X(B_k(t))\neq X(B_{k+1}(t))$ and $X(t)$ exists. As
$|C'(B_k(t))|$ is odd for all $k$ and $t$ the limit must be either
1 or $-1$.
\end{proof}

%

We are now ready to prove Lemma \ref{prop:main}.
\begin{proof}[Proof of Lemma \ref{prop:main}]
For every $t$ and $s$ $B_k(t)=B_k(s)$ for all $k$ sufficiently
large. Therefore $X(t)=X(s)$ for all $s$ and $t$. (See Figure
\ref{fig:arrows} at page \pageref{fig:arrows}.) Since $\mu$ is
ergodic and $X(0)$ is shift invariant, $X(0)$ is a $\mu$ almost
sure constant, which we denote by $X(\mu)$. Let
\begin{equation*}
\tilde \mu(A)=\mu(\{(x,y):(-x,y) \in A \}).
\end{equation*}
Then, by Lemma \ref{lem:symmetry}, $\tilde \mu$ is also a Gibbs
measure for $g$. On the other hand, $X(\tilde \mu)=-X(\mu)\neq
X(\mu)$, and therefore $\mu\neq\tilde\mu$.
\end{proof}

\section{An open problem}\label{sec:oppr}

Our construction is not monotone - by changing a $y$ value from a
$-1$ to a $+1$ we may change the set $S$ at which we look, and
then reduce the probability that the function outputs $+1$ at the
$x$ coordinate. The construction of Bramson and Kalikow, on the
other hand, is monotone. Therefore we ask the following question:
\begin{question}
Is there a value of $p$ and a continuous monotone regular \spec
$g$ such that $\var(g)\in\ell^p$ and $g$ admits multiple Gibbs
measures?
\end{question}


\noindent {\bf Acknowledgment:} Our research on this problem
started while N.B. and C.H. visited V.S. in IMPA, Rio de Janeiro.
We thank Yuval Peres and Maria Eulalia Vares for useful
discussions. N.B. and C.H. are grateful for the hospitality and
financial support during their stay at IMPA. The work was
partially supported by grants of Faperj and CNPq.

%

\end{document}